\documentclass[reqno,11pt,twoside]{article}
\usepackage{amssymb,latexsym,amsmath,amsthm,epsf,graphicx,graphics,psfrag}
\theoremstyle{plain}
\newtheorem{theorem}{Theorem}
\newtheorem{lemma}[theorem]{Lemma}
\newtheorem{proposition}[theorem]{Proposition}
\newtheorem{corollary}[theorem]{Corollary}

\newtheorem*{conjecture*}{Conjecture}
\newtheorem*{proposition*}{Proposition}
\newtheorem*{no label}{}

\theoremstyle{definition}
\newtheorem{definition}[theorem]{Definition}
\newtheorem{example}[theorem]{Example}

\theoremstyle{remark}
\newtheorem*{obs}{Remark}

\numberwithin{theorem}{section}

\def\ring#1{\mathring#1}
\def\cal#1{\mathcal{#1}}

\def\irred{irreducible}
\def\hndl{handlebody}
\def\hndls{handlebodies}
\def\compression{compression body}

\def\figurefontsize{9}

\title{Examples of Irreducible Automorphisms of Handlebodies}
\author{Leonardo Navarro Carvalho\footnote{Partially supported
        by CNPq--Brazil, CAPES--Brazil, FAPESP--S\~{a}o Paulo,
        Brazil, 03/06914-2.}\\
        {\small leonardo@ime.unicamp.br}}
\date{January 2004}
\begin{document}

\maketitle

\begin{abstract}
Automorphisms of handlebodies arise naturally in the
classification of automorphisms of three-manifolds.  Among
automorphisms of handlebodies, there are certain automorphisms 
called  \irred\ (or generic), which are analogues of pseudo-Anosov
automorphisms of surfaces.  We show that \irred\ automorphisms of
handlebodies exist and develop methods for constructing a range of
examples.
\end{abstract}


\section{Introduction.}\label{S:intro}

\subsection{Some history and background}
The classification of {\em automorphisms} (i.e.
self-homeomorphisms) of a manifold, up to isotopy, is a natural
and important problem. Nielsen addressed the case where the
manifold is a compact and connected surface and his results were
later substantially improved by Thurston (see
\cite{Nielsen:SurfacesI,Nielsen:SurfacesII,Nielsen:SurfacesIII,WPT:SurfaceAutos,MHWT:Nielsen}).
We briefly state their main result: An automorphism of a surface
is, up to isotopy, either {\em periodic} (i.e., has finite order),
{\em reducible} (i.e., preserves an essential codimension-1
submanifold) or {\em pseudo-Anosov}.  We refer the reader to any
of
\cite{FLP:TravauxThurston,MHWT:Nielsen,WPT:SurfaceAutos,ACSB:SurfaceAutos}
for details --- including the definition of a pseudo-Anosov
automorphism. The Nielsen-Thurston theory also shows that the
reducible case may --- as expected --- be reduced to the other
two. Since periodic automorphisms are relatively easy to
understand, the remaining {\em irreducible} case --- pseudo-Anosov
--- is the most interesting and rich one.

Indeed, pseudo-Anosov automorphisms of surfaces are the subject of
intense and wide research (see \cite{WPT:SurfaceAutos}). We
mention two works on the natural problem of building examples of
such automorphisms: Penner provides a generating method
\cite{RP:88} and a testing algorithm is developed in
\cite{BH:Surfaces}.

In \cite{UO:Autos}, Oertel undertakes a similar classification
project for a certain class of three-dimensional manifolds.
Suppose that a three-dimensional manifold $M$ is compact,
connected, orientable and {\em irreducible} (i.e., every embedded
sphere bounds a ball). Assume further that $\partial
M\neq\emptyset$. By use of canonical decompositions of $M$
due to Bonahon (determined by his {\em characteristic compression
body} \cite{FB:CompressionBody}) and Jaco, Shalen and Johanson
(the {\em JSJ-decomposition}
\cite{WJPS:Characteristic,KJ:Characteristic}), the study of
automorphisms of $M$ is reduced to the study of automorphisms of
{\em compression bodies} and {\em handlebodies} (see
\cite{UO:Autos}). We define these types of manifolds:

\begin{definition}
A {\em handlebody} $H$ is an orientable and connected
three-manifold obtained from a three-dimensional ball by attaching
a certain finite number $g$ of 1-handles. The integer $g$ is the
{\em genus} of $H$. It should be clear that $\pi_1(H)$ is
isomorphic to the free group $F_g$ on $g$ generators.

A {\em compression body} is a pair $(Q,F)$, where $Q$ is a
three-manifold obtained from a compact surface $F$ (not
necessarily connected) in the following way: consider the disjoint
union of $F\times I$ with the disjoint union of finitely many
balls $B$ and add 1-handles to $(F\times\{1\})\cup\partial B$,
obtaining $Q$. We allow empty or non-empty $\partial F$, but $F$
cannot have sphere components. We identify $F$ with
$F\times\{0\}\subseteq Q$, which is called the {\em interior
boundary} of $(Q,F)$, denoted by $\partial_i Q$. The {\em exterior
boundary} $\partial_e Q$ of $Q$ is the closure $\overline{\partial
Q-\partial_i Q}$. If $Q$ is homeomorphic to the disjoint union of
$F\times I$ with balls then $(Q,F)$ is said to be {\em trivial}.

We may abuse notation and refer to $Q$ as a compression body.
\end{definition}

The role of the disjoint union of balls $B$ in the definition of
\compression\ above is a two fold one. It makes some operations of
attaching 1-handles to be trivial. For instance, a 1-handle may
connect $F\times I$ with a ball or connect two distinct balls. On
the other hand, under this definition, a handlebody may be
regarded as a connected \compression\ whose interior boundary
$\partial_i Q=F$ is empty. Indeed, these two types of manifolds
are quite similar \cite{FB:CompressionBody}, as is the study of
their automorphisms \cite{UO:Autos}. This research will focus on
the case of handlebodies.

The following definition is due to Oertel:

\begin{definition}\label{D:red,irred}
An automorphism $f\colon H\to H$ of a handlebody $H$ is said {\em
reducible} if any of the following holds:
\begin{enumerate}
   \item there exists an $f$-invariant (up to isotopy) non-trivial
   compression body $(Q,F)$ with $Q\subseteq H$, $\partial_e
   Q\subseteq \partial H$ and $F=\partial_i Q\neq\emptyset$ not
   containing $\partial$-parallel disc components,

   \item there exists an $f$-invariant (up to isotopy) collection of
   pairwise disjoint, incompressible, non-$\partial$-parallel and
   properly embedded annuli, or

   \item $H$ admits an $f$-invariant (up to isotopy) $I$-bundle
   structure.
\end{enumerate}

The automorphism $f$ is said {\em \irred} (or {\em generic}, as in
\cite{UO:Autos}) if both of the following conditions hold:
\begin{enumerate}
    \item $\partial f=f|_{\partial H}$ is pseudo-Anosov, and

    \item there exists no {\em closed reducing surface} $F$:
    a closed reducing surface is a surface $F\neq\emptyset$ which
    is the interior boundary $\partial_i Q$ of a non-trivial
    compression body $(Q,F)$ such that $Q\subseteq H$, $(Q,F)$
    is $f$-invariant (up to isotopy) and $\partial_e Q=\partial H$.
\end{enumerate}
\end{definition}

An obvious remark is that this definition of \irred\ automorphism
excludes the periodic case.

\begin{theorem}[Oertel, \cite{UO:Autos}]\label{T:oertel classification}
An automorphism of a handlebody is either:
\begin{enumerate}
    \item periodic,

    \item reducible, or

    \item \irred.
\end{enumerate}
\end{theorem}

We note that the theorem above is not entirely obvious. For
example, one must show that if an automorphism $f:H\to H$ of a
handlebody does not restrict to a pseudo-Anosov $\partial f$ on
$\partial H$, then $f$ is actually reducible according to
Definition~\ref{D:red,irred}, or periodic.

Our interest is precisely in the \irred\ case, which is in many
ways analogous to the pseudo-Anosov case for surfaces (an
important similarity is related to the existence of certain
invariant projective measured laminations \cite{UO:Autos}).

\subsection{This research}
In this article we address the problem of constructing examples of
\irred\ automorphisms of handlebodies. We note that no examples of
such automorphisms were known before this research was undertaken.
The examples will give some indication that the theory of \irred\
automorphisms of handlebodies is even richer than the theory of
pseudo-Anosov automorphisms of surfaces.

Our results will give sufficient conditions for an automorphism to
be irreducible. These conditions will either be constructible or
verifiable, so they can (and will) be used to generate actual examples.

In Section \ref{S:first example} we build an example of an \irred\
automorphism of a genus two handlebody (see Example~\ref{E:first
example} and Proposition~\ref{P:reducao bordo}).

In Section \ref{S:method} we generalize the construction of that
first example and develop a method for generating a larger class
of \irred\ automorphisms. This is done in theorems \ref{T:method}
and \ref{T:method2}. Their statements depend on some rather
technical constructions, unsuited for this introduction.

The final section concerns closed reducing surfaces. Under some
hypotheses, we will find bounds (upper and lower) for the Euler
characteristic of a possible closed reducing surface. As a
corollary we shall show:
\begin{no label}[Corollary~\ref{C:genus2}]
Let $f\colon H\to H$ restrict to $\partial H$ as a pseudo-Anosov
automorphism. If the genus of $H$ is two then $f$ is \irred.
\end{no label}
One can use this as a tool to build examples of \irred\
automorphisms. A consequence will be that the complexity of an
\irred\ automorphism of a handlebody cannot be extracted from the
homomorphism induced in the fundamental group
--- unlike the analogous situation in dimension two.
Example~\ref{E:genus2} will describe an \irred\ automorphism of a
handlebody $H$ whose induced automorphism in $\pi_1(H)\to\pi_1(H)$
is the identity.

The following two theorems will serve us as important tools. We
refer the reader to \cite{RP:88} and \cite{BH:Tracks} respectively
for details and precise definitions.

\begin{theorem}[Penner, \cite{RP:88}]\label{T:penner}
Let $\mathcal{C}$, $\mathcal{D}$ be two systems of closed curves
in an orientable surface $S$ with $\chi(S)<0$. Assume that
$\cal{C}$ and $\cal{D}$ intersect efficiently, do not have
parallel components and fill $S$.  Let $f\colon S\to S$ be a
composition of Dehn twists: right twists along curves of
$\mathcal{C}$ and left twist along curves of $\mathcal{D}$. If a
twist along each curve appears at least once in the composition,
then $f$ is isotopic to a pseudo-Anosov automorphism of $S$.
\end{theorem}

\begin{theorem}\label{T:pA irreducible}
Let $S$ be a compact surface with $\chi(S)<0$ and precisely one
boundary component. An automorphism $f\colon S\to S$ is
pseudo-Anosov if and only if $f_*^n$ is irreducible for all $n>0$.
\end{theorem}

As a final introductory remark, we observe that an ideal
classification of automorphisms of handlebodies should identify in
each isotopy class a representative which is ``best" in some
sense. Considering the classification of Theorem~\ref{T:oertel
classification}, this has been done for periodic and many
reducible automorphisms \cite{FB:CompressionBody,UO:Autos}. The
problem of finding a best representative of an irreducible
automorphism is addressed in \cite{UO:Autos,LCarvalho:thesis} but
not yet solved.

We will adopt the following notation: given a topological space
$A$ (typically a manifold or sub-manifold), $\overline{A}$ will
denote its topological closure, $\mathring{A}$ its interior and
$|A|$ its number of connected components. If $M$ is a manifold and
$S\subseteq M$ a compact codimension 1 submanifold, we can ``cut
$M$ open along $S$'' obtaining $M_S$. More precisely, a Riemannian
metric in $M$ determines a path-metric in $M-S$, in which the
distance between two points is the infimum of the lengths of paths
in $M-S$ connecting them. We let $M_S$ to be the completion of
$M-S$ with this metric.

I thank Ulrich Oertel for many enlightening meetings and helpful
suggestions in his role as dissertation advisor, and also for
laying the foundations on which the research in this paper is
built.


\section{An example.}\label{S:first example}

We show that \irred\ automorphisms of \hndls\ exist by presenting
an example.

Let $H$ be a genus 2 \hndl. We will describe an automorphism of
$H$ as a composition of Dehn twists along two annuli and a disc.
We shall prove that it is \irred\ by showing that its restriction
to $\partial H$ is pseudo-Anosov and that, for an algebraic
reason, there can be no closed reducing surface.

%
%

\begin{example}\label{E:first example}
We start with a pseudo-Anosov automorphism $\varphi\colon S\to S$
of the once punctured torus $S$. Such a $\varphi$ will be defined
as a composition of Dehn twists along two curves.

We will represent $S$ as a cross, after identifying pairs of
opposite sides as shown in  Figure~\ref{F:toro}.

\begin{figure}[h]
\centering
\psfrag{alpha0}{\fontsize{\figurefontsize}{12}$\alpha_0$}
\psfrag{alpha1}{\fontsize{\figurefontsize}{12}$\alpha_1$}
\includegraphics[scale=0.28]{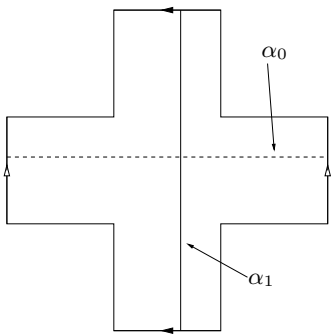}
\caption{\small The oriented surface $S$ and the curves
$\alpha_0$, $\alpha_1$.}\label{F:toro}
\end{figure}

Let $\alpha_0$, $\alpha_1$ be simple closed curves as in the
figure. It is easy to verify that the systems $\cal
C=\{\alpha_0\}$ e $\cal D=\{\alpha_1\}$ satisfy the hypothesis of
Theorem~\ref{T:penner} (Penner). Let $T_0^-$ be the left Dehn
twist along $\alpha_0$ and $T_1^+$ the right twist along
$\alpha_1$. We define:
$$
\varphi=T_1^+\circ T_0^-.
$$
By Theorem~\ref{T:penner}, $\varphi$ is pseudo-Anosov. Then, by
Theorem~\ref{T:pA irreducible}, any positive power $\varphi_*^n$,
of the induced homomorphism $\varphi_*\colon\pi_1(S)\to\pi_1(S)$
is irreducible. We note this fact for future use.

We now consider the \hndl\ $H=S\times I$ and lift $\varphi$ to
$H$, obtaining $\phi\colon H\to H$, a composition of twists along
the annuli $A_0=\alpha_0\times I$, $A_1=\alpha_1\times I$ as in
Figure~\ref{F:aneis}.

\begin{obs}
For future use, we will think of the picture as being looked at
``from above''. More precisely, we orient $H$ in such a way that
the induced orientation in $S\times\{1\}$ coincides with the one
inherited naturally from $S$.
\end{obs}

Identifying $\pi_1(H)$ with $\pi_1(S)$ we have $\phi_*=\varphi_*$.

\begin{figure}[h]
\centering \psfrag{Delta}{\fontsize{\figurefontsize}{12}$\Delta$}
\psfrag{D}{\fontsize{\figurefontsize}{12}$\mathcal{C}$}
\psfrag{E}{\fontsize{\figurefontsize}{12}$\mathcal{D}$}
\psfrag{A0}{\fontsize{\figurefontsize}{12}$A_0$}
\psfrag{A1}{\fontsize{\figurefontsize}{12}$A_1$}
\psfrag{H}{\fontsize{\figurefontsize}{12}$H$}
\includegraphics[scale=0.3]{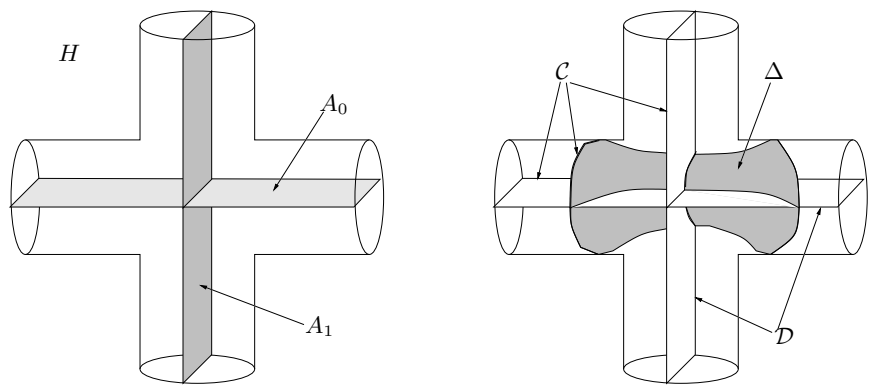} \caption{\small The
automorphism $f$ is defined as a composition of Dehn twists along
the annuli $A_0$, $A_1$ and the disc $\Delta$.} \label{F:aneis}
\end{figure}

Finally, we will obtain the desired \irred\ automorphism $f\colon
H\to H$ by composing $\phi$ with a twist along a disc $\Delta$,
shown in Figure~\ref{F:aneis}.

Let $T_\Delta^+$ be the right Dehn twist along $\Delta$. We
define:
$$
f=T_\Delta^+\circ\phi.
$$

\begin{proposition}\label{P:reducao bordo}
The automorphism $f\colon H\to H$ is \irred.
\end{proposition}

Part of the proof will be done in the following general lemma.

\begin{lemma}\label{L:reducible pi_1}
Let $g\colon H\to H$ be an automorphism of a \hndl\ $H$ such that
$\partial g$ is pseudo-Anosov. If $g$ is reducible then, for some
$n\in\mathbb{N}$, $g_*^n\colon\pi_1(H)\to\pi_1(H)$ is reducible.
\end{lemma}
\begin{proof}
Let $Q$ be a \compression\ invariant under $g$. Let
$F\subseteq\partial_iQ$ be a component of the closed reducing
surface and $J\subseteq\ring H$ the \hndl\ bounded by $F$.
Choosing a base point in $J$ and omitting the obvious inclusion
homomorphisms we claim that
$$
\pi_1(H)=\pi_1(J)*G,
$$
where $G$ is not trivial. To see this, first consider the
connected and nontrivial \compression\ $Q'=\overline{H-J}$, whose
boundary decomposes as $\partial_i Q'=F$ and $\partial_e
Q'=\partial H$. The \compression\ structure of $Q'$ gives it as a
product $F\times I$ to which 1-handles are attached. Regarding
$F\times I\subseteq Q'\subseteq H$, we see that the \hndl\
$J'=(F\times I)\cup J$ deformation retracts to $J$ (so
$\pi_1(J')=\pi_1(J)$ through inclusion). But the \compression\
structure of $Q'$ gives $H$ as $J'$ with 1-handles attached to
$\partial J'$. Since $\partial J'$ is connected, we can moreover
assume that these 1-handles are attached to a disc in $\partial
J'$, which gives $\pi_1(H)=\pi_1(J')*G=\pi_1(J)*G$, where $G$ is a
free group (whose rank equals the number of 1-handles of $Q'$).
Since $Q'$ is not trivial, $G$ is not trivial, proving the claim.
Therefore $\pi_1(J)$ is a proper free factor of $\pi_1(H)$.

Let $g^n$ be the first power of $g$ preserving $J$. Isotoping $g$
we assume moreover that the base point is fixed by $g^n$. From
$$
g^n(J)=J
$$
follows that $g_*^n(\pi_1(J))$ is conjugate to $\pi_1(J)$, hence
the class of $g_*^n$ in $Out\big(\pi_1(H)\big)$ is reducible.
\end{proof}

\begin{proof}[Proof of Proposition~\ref{P:reducao bordo}]
We need to prove that $\partial f=f|_{\partial H}$ is
pseudo-Anosov and that $f$ does not admit closed reducing
surfaces.

We start by verifying that $\partial f$ is pseudo-Anosov. It is
given as composition of Dehn twists: left twists along curves of
$$
\mathcal{C}=\left\{\,(\alpha_0\times\{1\})\,,\,(\alpha_1\times\{0\})\,\right\},
$$
(see Figure~\ref{F:aneis}) and right twists along curves of
$$
\mathcal{D}=\left\{\,(\alpha_0\times\{0\})\,,\,(\alpha_1\times\{1\})\,,\,\partial\Delta\,\right\}.
$$
We now note that $\mathcal{C}$, $\mathcal{D}$ satisfy the
hypotheses of Theorem~\ref{T:penner}, hence $\partial f$ is
pseudo-Anosov.

We prove by contradiction that $f$ admits no closed reducing
surface.  Suppose there is a closed reducing surface.  By
Lemma~\ref{L:reducible pi_1}, there exists $n$ such that $f_*^n$
is reducible. But $f=(T_\Delta^+)\circ\phi$ and the twist
$(T_\Delta^+)$ (along a disc) induces the identity in $\pi_1(H)$.
Therefore, recalling that $\pi_1(H)$ is identified with
$\pi_1(S)$, we have that $f_*^n=\phi_*^n=\varphi_*^n$, which was
seen before to be irreducible for any $n$, a contradiction.

Therefore $f$ is \irred.
\end{proof}
\end{example}


\section[A method]{A method for generating \irred\ automorphisms.}
\label{S:method}

The construction of Example~\ref{E:first example} may be
generalized to provide a method for generating a larger class of
\irred\ automorphisms of \hndls\ (Theorems \ref{T:method} and
\ref{T:method2}). This method partially solves a problem proposed
in \cite{UO:Autos}.

\begin{definition}\label{D:penner pair}
We say that a pair $(\cal{C},\cal{D})$ of curve systems in a
compact, connected and orientable surface $S$ with $\chi(S)<0$ is
a {\em Penner pair in $S$} if $\cal C$, $\cal D$ satisfy the
hypotheses of Penner's Theorem~\ref{T:penner} i.e.,
\begin{enumerate}
    \item each $\cal C$, $\cal D$ is a finite collection of simple,
    closed and pairwise disjoint essential curves,

    \item $\cal{C}$ and $\cal{D}$ intersect efficiently, do not
    have parallel components and {\em fill} $S$
    (i.e., the components of $S-(\cal{C}\cup\cal{D})$ are either
    contractible or deformation retract to $\partial S$).
\end{enumerate}
Suppose that $(\cal C,\cal D)$ is a Penner pair. An automorphism
$\varphi$ of $S$ obtained from $\cal C$, $\cal D$ as in
Theorem~\ref{T:penner} is called a {\em Penner automorphism
subordinate to $(\cal C,\cal D)$}.

If $\partial S\neq\emptyset$ then a properly embedded and
essential arc $\theta$ is called {\em dual to $ (\cal C,\cal D)$}
if $\theta$ intersects $\cal C\cup\cal D$ transversely and in
exactly one point $p\notin\cal C\cap\cal D$.
\end{definition}

\begin{obs}
Although not all Penner pairs admit dual arcs it is easy to
construct pairs that do: such a pair $(\cal C,\cal D)$ in $S$ has
the property that there are two adjacent components (not
necessarily distinct) of $S-(\cal C\cup\cal D)$ each containing
some component of $\partial S$. If a pair does not have this
property then we can remove discs from $S$ and introduce dual
arcs.
\end{obs}

We constructed the \irred\ automorphism in Example \ref{E:first
example} by lifting a pseudo-Anosov automorphism of a surface to a
product and composing it with a twist on a disc. The general
method will be similar. Our interest in dual arcs is that we can
use them to construct discs that will yield \irred\ automorphisms.

Throughout this section we fix a compact, connected and oriented
surface $S$ with $\partial S\neq\emptyset$ and define $H=S\times
I$, which is a handlebody. We identify $S$ with
$S\times\{1\}\subseteq H$, inducing orientation in $H$.

Given a Penner pair $(\cal C,\cal D)$ in $S$ and a dual arc
$\theta$, we build a disc $\Delta_\theta$ in $H$ in the following
way. Let $\gamma$ be the curve of $(\cal C,\cal D)$ that $\theta$
intersects and assume without loss of generality that
$\gamma\subseteq\cal{C}$. Let $D=\theta\times I\subseteq H$. Then
$\partial D$ intersects $\gamma_1=\gamma\times\{1\}$ in a point.
Now let $\Delta_\theta$ be the {\em band sum} of $D$ with itself
along $\gamma_1$.  This means that $\Delta_\theta$ is obtained from $D$ and
$\gamma_1$ by the following construction: consider a regular
neighborhood $N=N(D\cup\gamma_1)$. Then
$\Delta_\theta=\overline{\partial N-\partial H}$ is a properly
embedded disc.

\begin{theorem}\label{T:method}
Suppose that $\partial S\neq\emptyset$ has exactly one component.
Let $(\cal C,\cal D)$ be a Penner pair in $S$ with dual arc
$\theta$ and $\varphi\colon S\to S$ a Penner automorphism
subordinate to $(\cal C,\cal D)$. Let $\hat\varphi\colon H\to H$
be the lift of $\varphi$ to the product $H=S\times I$ and
$\Delta_\theta\subseteq H$ the disc constructed from the arc
$\theta$ as above. Then there exists a Dehn twist
$T_{\Delta_\theta}\colon H\to H$ along $\Delta_\theta$ such that
the composition
$$
\hat\varphi\circ T_{\Delta_\theta} \colon H\to H
$$
is an \irred\ automorphism of $H$.
\end{theorem}

The key to the proof is the verification that $\cal{C}$, $\cal{D}$ and
$\partial\Delta_\theta$ induce a Penner pair in $\partial H$.

\begin{lemma}\label{L:twisting disc}
Let $S$, $(\mathcal{C},\mathcal{D})$, $\theta$, $H=S\times I$ and
$\Delta_\theta$ be as in the statement of Theorem~\ref{T:method}.
Let $\cal C_i=\cal C\times\{i\}\subseteq S_i=S\times\{i\}$ and
$\cal D_i=\cal D\times\{i\}\subseteq S_i=S\times\{i\}$, defining
$\cal C_0$, $\cal D_0\subseteq S_0$ and $\cal C_1$, $\cal
D_1\subseteq S_1$. Under these conditions the following system of
curves in $\partial H$:
\begin{align}
\mathcal{Q} =& \cal D_0 \cup \cal C_1 \cup\{\,\partial\Delta_\theta\,\},\notag\\
\mathcal{R} =& \cal C_0 \cup \cal D_1 \notag,
\end{align}
determine a Penner pair $(\mathcal{Q},\mathcal{R})$ in $\partial
H$.
\end{lemma}
\begin{proof}
We start by making the obvious remarks that $\cal C_0$, $\cal
D_0$, $\cal C_1$, $\cal D_1\subseteq\partial H$ and
$\cal{C}_0\cap\cal{D}_1=\emptyset$,
$\cal{D}_0\cap\cal{C}_1=\emptyset$. Recall we are assuming that
$\theta\cap(\cal{C}\cup\cal{D})\subseteq\gamma\subseteq\cal{C}$.
We verify that:
\begin{itemize}

    \item $\partial\Delta_\theta\cap\cal{D}_0=\emptyset,$ because
    $\left(\theta\times\{0\}\right)\cap\cal{D}_0=\emptyset$ and
    $\partial\Delta_\theta\cap S_0$ consists of two arcs parallel to
    $\theta\times\{0\}$,

    \item $\partial\Delta_\theta\cap\cal{C}_1=\emptyset,$ because
    $\partial\Delta_\theta\cap\gamma_1=\emptyset$ by construction.
\end{itemize}
Therefore each $\mathcal{Q} =\cal D_0 \cup \cal
C_1\cup\{\partial\Delta\}$ and $\mathcal{R} = \cal C_0 \cup\cal
D_1$ is a system of simple closed curves essential in $\partial
H$. To conclude that $(\mathcal{Q},\mathcal{R})$ is indeed a
Penner pair we just need to verify that
$\mathcal{Q}\cup\mathcal{R}$ fills $\partial H$.

A component of $S-(\cal C\cup\cal D)$ either is a disc or an
annulus that retracts to $\partial S$. Therefore a component of
$\partial H-(\cal C_0\cup\cal D_0\cup \cal C_1\cup\cal D_1)$
either is a disc, or an annulus $A$ (that retracts to $\partial
S\times I$). But $A\cap\partial\Delta_{\theta}$ is a union of four
arcs essential in $A$, hence each component of $\partial
H-(\mathcal{Q}\cup\mathcal{R})$ is a disc, showing that
$\mathcal{Q}\cup\mathcal{R}$ fills $\partial H$, completing the
proof.
\end{proof}

Instead of proving Theorem \ref{T:method} we will prove the more
general result below, which clearly implies the other. We note
that twists on curves of $\cal C$, $\cal D$ in $S$ lift to twists
along annuli in $H$. We call these systems of annuli $\hat{\cal
C}$, $\hat{\cal D}$ respectively.  It thus makes sense to refer to
directions of the twists along these vertical annuli (recall that
$H$ has orientation induced by $S\times\{1\}\subseteq H$).

\begin{theorem}\label{T:method2}
Let $(\mathcal{C},\mathcal{D})$, $S$, $\theta$, $H$ and
$\Delta_\theta$ be as in Theorem~\ref{T:method}. Let $f$ be a
composition $f\colon H\to H$ of twists along the annuli of
$\hat{\cal C}$, $\hat{\cal D}$ and the disc $\Delta_{\theta}$: in
one direction along the annuli in $\hat{\cal D}$ and in the
opposite direction along the annuli in $\hat{\cal C}$ and the disc
$\Delta_{\theta}$. If each of these twists appear in the
composition at least once $f$ is \irred.
\end{theorem}
\begin{proof}
We will show initially that $f_*^n\colon\pi_1(H)\to\pi_1(H)$ is an
irreducible automorphism of a free group for any $n\geq 0$ (hence
there can be no closed reducing surface by Lemma~\ref{L:reducible
pi_1}) and then that $\partial f=f|_{\partial H}$ is
pseudo-Anosov, thus completing the proof that $f$ is \irred.

We first identify $\pi_1(H)$ with $\pi_1(S)$, identifying $S$ with
$S\times\{1\}\subseteq H$. Let $T_{\Delta_{\theta}}$ be a twist
along $\Delta_{\theta}$. Since
$(T_{\Delta_\theta})_*\colon\pi_1(H)\to\pi_1(H)$ is the identity
($\Delta_\theta$ is a disc) the hypotheses on $f$ imply that
$f_*=\varphi_*$ for some Penner automorphism $\varphi\colon S\to
S$ subordinate to $(\cal C,\cal D)$. Penner automorphisms are
pseudo-Anosov so, given that $\partial S$ has a single component,
it follows from Theorem~\ref{T:pA irreducible} that $\varphi_*^n$
is an irreducible automorphism of $\pi_1(S)$ for any $n\geq 0$.
Therefore $f_*^n\colon\pi_1(H)\to\pi_1(H)$ is irreducible, proving
that $f$ does not admit closed reducing surfaces
(Lemma~\ref{L:reducible pi_1}).

We now prove that $\partial f$ is pseudo-Anosov. Let
$(\mathcal{Q},\mathcal{R})$ be as in Lemma~\ref{L:twisting disc},
therefore a Penner pair. By construction the twists that compose
$f$ restrict to $\partial H$ as twists along curves of $\cal{Q}$
or $\cal{R}$. It is then straightforward to verify that $\partial
f$ is a Penner automorphism subordinate to $(\cal{Q},\cal{R})$,
hence pseudo-Anosov, completing the proof that $f$ is irreducible.
\end{proof}

\begin{example}\label{E:method}
Consider $S$ a genus 2 surface minus a disc, represented in
Figure~\ref{F:twisting curves} as an octagon whose sides are
identified according to the arrows.

\begin{figure}[ht]
\centering \psfrag{aa}{\fontsize{\figurefontsize}{12}$\alpha$}
\psfrag{bb}{\fontsize{\figurefontsize}{12}$\beta$}
\psfrag{cc}{\fontsize{\figurefontsize}{12}$\gamma$}
\psfrag{dd}{\fontsize{\figurefontsize}{12}$\delta$}
\psfrag{tt}{\fontsize{\figurefontsize}{12}$\theta$}
\psfrag{a}{\fontsize{\figurefontsize}{12}$A$}
\psfrag{b}{\fontsize{\figurefontsize}{12}$B$}
\psfrag{c}{\fontsize{\figurefontsize}{12}$C$}
\psfrag{d}{\fontsize{\figurefontsize}{12}$D$}
\includegraphics[scale=0.35]{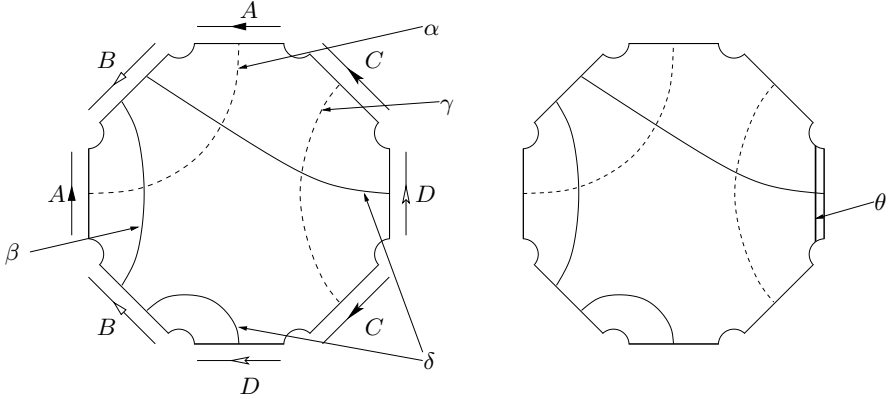} \caption{\small
A Penner pair in $S$, with dual arc $\theta$.} \label{F:twisting
curves}
\end{figure}

In the picture there are represented four further curves:
$\alpha$, $\beta$, $\gamma$ and $\delta$. Defining
\begin{align}
\cal C=&\{\,\beta,\delta\,\},\notag\\
\cal D=&\{\,\alpha,\gamma\,\},\notag
\end{align}
it is easy to check that $(\cal C,\cal D)$ is a Penner pair in
$S$. The automorphism $\varphi\colon S\to S$ defined by
$$
\varphi=T_\beta^-\circ T_\delta^-\circ T_\alpha^+\circ T_\gamma^+
$$
is, therefore, a Penner automorphism subordinate to the pair
$(\cal C,\cal D)$.

The pair $(\cal C,\cal D)$ admits dual arcs. The picture shows
one, labelled as $\theta$. We consider the corresponding disc
$\Delta_{\theta}$. Figure~\ref{F:twists in product} shows
$S_0=S\times\{0\}$, $S_1=S\times\{1\}\subseteq\partial H$ and how
$\partial\Delta_{\theta}$ intersects them.

\begin{figure}[ht]
\centering
\psfrag{SX0}{\fontsize{\figurefontsize}{12}$S\times\{0\}$}
\psfrag{SX1}{\fontsize{\figurefontsize}{12}$S\times\{1\}$}
\psfrag{pD}{\fontsize{\figurefontsize}{12}$\partial\Delta_{\theta}$}
\includegraphics[scale=0.35]{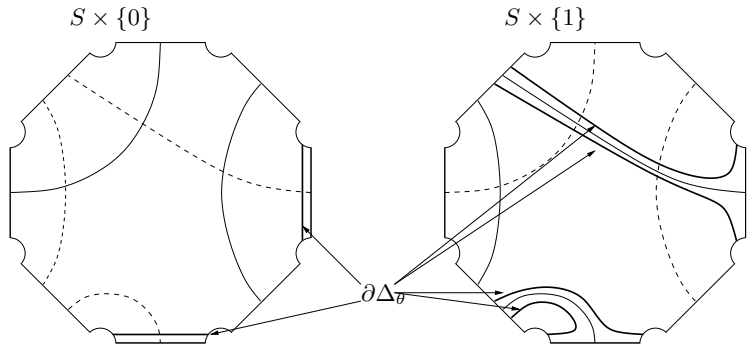}
\caption{\small The curve $\partial\Delta_{\theta}$ in $\partial
H$.} \label{F:twists in product}
\end{figure}

Figure~\ref{F:twists in product}  shows the pair
$(\mathcal{Q},\mathcal{R})$ obtained by Lemma~\ref{L:twisting
disc} as well: $\mathcal{Q}$ consists on the solid curves,
including $\partial\Delta_{\theta}$, while the dotted curves form
$\mathcal{R}$.

Theorem~\ref{T:method} assures that, if $\hat\varphi\colon H\to H$
is the lift of $\varphi$ to $H$, then
$$
\hat\varphi\circ T_{\Delta_{\theta}} \colon H\to H
$$
is an \irred\ automorphism for a certain twist
$T_{\Delta_{\theta}}$ along $\Delta_{\theta}$.
\end{example}

\section{Genus of a closed reducing surface.}\label{S:genus2}

Consider an \irred\ automorphism $f$ of a handlebody $H$. By
definition, the restriction $\partial f=f|_{\partial H}$ is
pseudo-Anosov. We prove in this section that, for genus 2
handlebodies, the converse is also true: if $\partial f$ is
pseudo-Anosov $f$ is \irred\ (Corollary~\ref{C:genus2}). Therefore
methods for generating pseudo-Anosov automorphisms may be used for
generating irreducible automorphisms of handlebodies (e.g.,
Example~\ref{E:genus2}). The nonexistence of a closed reducing
surface in this case comes from a geometric argument: see
Theorem~\ref{T:reduction}, which determines bounds on the Euler
characteristic of a closed reducing surface. Therefore this
criterion does not depend on the automorphism induced on the
fundamental group and, for this reason, yields interesting
examples, even if the hypothesis on the genus looks restrictive.

The following lemma is elementary, so we leave the proof to the
reader.
\begin{lemma}\label{L:euler}
Let $(Q,F)$ be a connected and non-trivial \compression\ with
$F\neq\emptyset$. If $n$ is the smallest number of 1-handles of
$Q$ then
$$
\chi(\partial_eQ)=\chi(F)-2n.
$$
\end{lemma}

\begin{proposition}\label{P:1-handle}
Let $(Q,F)$ be a non-trivial \compression\ with only one
$1$-handle and $F\neq\emptyset$. Then the handle is unique up to
isotopy.
\end{proposition}

\begin{proof}
We start by noting that, up to isotopy, the $1$-handle and a {\em
dual disc} (the co-core of a $1$-handle) determine each other.

We fix a $1$-handle for $(Q,F)$ and a dual disc $D$. Consider
another choice of $1$-handle and pick a dual disc $D'$. We can
assume without loss of generality that $Q$ is connected: $Q$ has
just one $1$-handle so all but one of the components (the one
containing the initial $1$-handle) are either products or balls.
Since $\partial_e Q\cap L$ is incompressible in such a product (or
ball) component $L$, the component containing the initial
$1$-handle must contain both $D$ and $D'$, hence the other
$1$-handle as well.

We shall show that $D$ and $D'$ are isotopic. As mentioned
previously, this implies the proposition. The proof will be done
in two steps. The first step simplifies $D\cap D'$ through
standard ``cutting and pasting'' methods, obtaining disjoint $D$
and $D'$. In the second step we will show that $D$ and $D'$ must
be parallel, as desired.

We begin the first step by perturbing $D'$ so that $D\cap D'$
consists of closed curves and arcs. We can eliminate closed curves
from $D\cap D'$ by standard ``cut and paste'' arguments. Since $Q$
is irreducible this can be attained by isotopy of $D'$. After a
finite number of such operations we have that $D\cap D'$ contains
no closed curves.

So assume that $D\cap D'$ consists of arcs. Here again we will
perform isotopies that will reduce $|D\cap D'|$. Recall that
$F=\partial_iQ$ and consider two cases: 1, that $F$ is
disconnected, and 2, that $F$ is connected.

\

1. $F$ is disconnected. In this case $D$ separates $Q$.  Let
$\alpha$ be an arc of $D\cap D'$ which is edgemost in $D'$,
cutting a half-disc $\Delta'$ from $D'$.

Let $Q'=\overline{Q-D}$ with product structure inherited from $Q$.
Since $D$ separates $Q$ it follows that $Q'$ has two components
(see Figure \ref{F:split1}).

\begin{figure}[h]
\centering \psfrag{Q}{\fontsize{\figurefontsize}{0}$Q$}
\psfrag{Q'}{\fontsize{\figurefontsize}{0}$Q'$}
\psfrag{D}{\fontsize{\figurefontsize}{0}$D$}
\psfrag{D0}{\fontsize{\figurefontsize}{0}$D_0$}
\psfrag{U0}{\fontsize{\figurefontsize}{0}$F_0$}
\psfrag{UXI}{\fontsize{\figurefontsize}{0}$F\times I$}
\psfrag{U}{\fontsize{\figurefontsize}{0}$F$}
\includegraphics[scale=0.4]{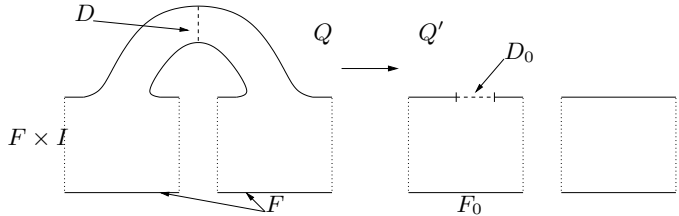} \caption{\small Cutting
$Q$ open along $D$ yields $Q'$.} \label{F:split1}
\end{figure}

Let $F_0\times I$ be the component of $Q'$ containing $\Delta'$.
There exists a disc $D_0\subseteq F_0\times\{1\}$ corresponding to
$D$ and, abusing notation, let $\alpha\subseteq D_0$ be defined by
$\alpha=\partial\Delta'\cap D_0$. Moreover, we have that
$\partial\Delta'\subseteq F_0\times\{1\}$. Since $F_0\times\{1\}$
is incompressible, $\partial\Delta'$ must bound a disc
$\Delta''\subseteq F_0\times \{1\}$. By the choice of $\Delta'$,
$\Delta''\cap D_0$ is a disc $\Delta'''$. Irreducibility of $Q'$
implies that the sphere $\Delta'\cup\Delta''$ bounds a ball.
Regluing $Q'$ to recover $Q$, such a ball determines another ball
$B$ in $Q$. We isotope $\Delta'$ through $B$ a little beyond
$\Delta'''$ to remove $\alpha$ from $D\cap D'$, thus reducing
$|D\cap D'|$.

A finite sequence of such isotopies yields $D\cap D'=\emptyset$.

\

2. $F$ is connected. Now $D$ does not separate $Q$. Consider again
an arc $\alpha\subseteq D\cap D'$ which is edge-most in $D'$,
bounding with an arc $\beta\subseteq\partial D'$ an edge-most disc
$\Delta'\subseteq D'$.

Let $Q'=\overline{Q-D}\simeq F\times I$ be the \compression\ $Q$
cut open along $D$ (Figure \ref{F:split2}).

\begin{figure}[h]
\centering \psfrag{Q}{\fontsize{\figurefontsize}{0}$Q$}
\psfrag{Q'}{\fontsize{\figurefontsize}{0}$Q'$}
\psfrag{D}{\fontsize{\figurefontsize}{0}$D$}
\psfrag{D+}{\fontsize{\figurefontsize}{0}$D^+$}
\psfrag{D-}{\fontsize{\figurefontsize}{0}$D^-$}
\psfrag{UXI}{\fontsize{\figurefontsize}{0}$F\times I$}
\psfrag{U}{\fontsize{\figurefontsize}{0}$F$}
\includegraphics[scale=0.4]{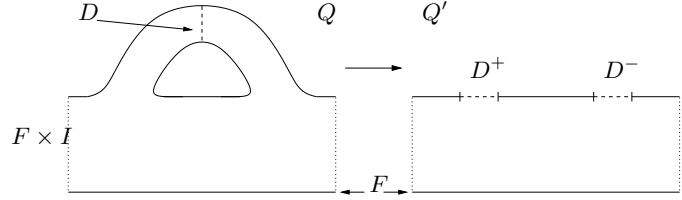} \caption{\small Cutting
$Q$ open along $D$ yields $Q'$.} \label{F:split2}
\end{figure}

In this case $F\times\{1\}$ contains two discs $D^+$, $D^-$
corresponding to $D$ and we suppose that $D^+$ corresponds to the
side of $D$ associated to $\Delta'$. There are two arcs,
$\alpha^+\subseteq D^+$ and $\alpha^-\subseteq D^-$ corresponding
to $\alpha$. Proceeding as in the previous case, we have that
$\partial\Delta'\cap D^+=\alpha^+$, $\partial\Delta'\cap
D^-=\emptyset$, $\partial\Delta'\subseteq F\times\{1\}$. Again
incompressibility gives us $\partial\Delta'$ bounding a disc
$\Delta''\subseteq F\times\{1\}$. There are two possibilities to
be considered depending on $\Delta''\cap D^-$:

(a) If $\Delta''\cap D^-=\emptyset$ (see Figure \ref{F:FX1}
    (a), recalling that $\alpha\cup\beta=\partial\Delta'=\partial\Delta''$)
    we proceed as in case 1, removing the arc
    $\alpha$ from $D\cap D'$, reducing $|D\cap D'|$.

(b) If $\Delta''\cap D^-\neq\emptyset$ (Figure \ref{F:FX1}
    (b)) the argument is more laborious. Since
    $\partial\Delta'\cap D'=\emptyset$, in this case we must have
    $\Delta''\cap D^-=D^-$.

\begin{figure}[h]
\centering \psfrag{pDelta}{\fontsize{\figurefontsize}{0}$\beta$}
\psfrag{Delta}{\fontsize{\figurefontsize}{0}$\Delta''$}
\psfrag{a}{\fontsize{\figurefontsize}{0}$\alpha$}
\psfrag{D+}{\fontsize{\figurefontsize}{0}$D^+$}
\psfrag{D-}{\fontsize{\figurefontsize}{0}$D^-$}
\psfrag{UX1}{\fontsize{\figurefontsize}{0}$F\times\{1\}$}
\includegraphics[scale=0.55]{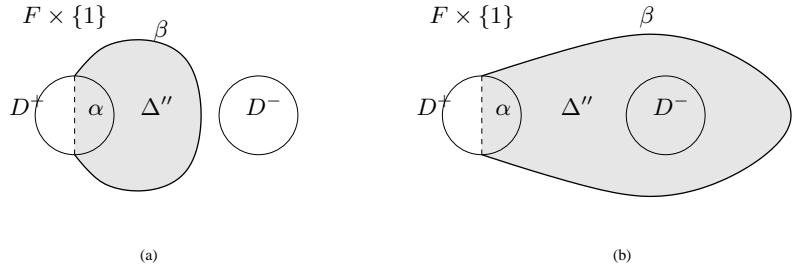} \caption{\small The disc
$\Delta''\subseteq F\times\{1\}$ parallel
    to $\Delta'$. The cases (a) and (b) depend on $\Delta''\cap
    D^-$.}
\label{F:FX1}
\end{figure}

Back to considering $D$, $D'\subseteq Q$, we note that $\partial
D'\cap D$ is finite. Therefore, the set
$\Gamma=\{\gamma_0,\gamma_1,\dots,\gamma_{n-1}\}$ of the closures
of the components of $\partial D'-D$ consists of finitely many
arcs.  We fix an arbitrary orientation for $\partial D'$ and use
it to induce a cyclic order in $\Gamma$. We can assume that the
indices respect this order and, since $\beta\in\Gamma$, we can
assume further that $\gamma_0=\beta$ (see Figure \ref{F:order}).

\begin{figure}[h]
\centering \psfrag{D'}{\fontsize{\figurefontsize}{0}$D'$}
\psfrag{b=g0}{\fontsize{\figurefontsize}{0}$\beta=\gamma_0$}
\psfrag{g1}{\fontsize{\figurefontsize}{0}$\gamma_1$}
\psfrag{g2}{\fontsize{\figurefontsize}{0}$\gamma_2$}
\psfrag{...}{\fontsize{\figurefontsize}{0}$\dots$}
\psfrag{gn-1}{\fontsize{\figurefontsize}{0}$\gamma_{n-1}$}
\includegraphics[scale=0.45]{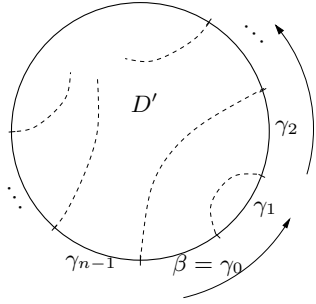} \caption{\small
Orientation of the arcs in $\partial D'$.} \label{F:order}
\end{figure}

In addition to the order in $\Gamma$, the orientation of $\partial
D'$ induces an orientation in each $\gamma_i$, giving an order
between the ends of each arc.

Working in $Q'$ again, assume that some arc $\gamma_i$ of $\Gamma$
has both extremes in $\partial D^-$. In this case $\gamma_i$, together
with an arc in $\partial D^-$, bounds a disc
$\Delta'''\subseteq\overline{(F\times\{1\}-(D^+\cup D^-))}$ for,
{\em in $Q'$}, $\gamma_i\cap D^+=\emptyset$,
$\gamma_i\cap\beta=\emptyset$ (see Figure \ref{F:delta3}).

\begin{figure}[h]
\centering \psfrag{pDelta}{\fontsize{\figurefontsize}{0}$\beta$}
\psfrag{Delta}{\fontsize{\figurefontsize}{0}$\Delta''$}
\psfrag{a}{\fontsize{\figurefontsize}{0}$\alpha$}
\psfrag{D+}{\fontsize{\figurefontsize}{0}$D^+$}
\psfrag{D-}{\fontsize{\figurefontsize}{0}$D^-$}
\psfrag{Delta3}{\fontsize{\figurefontsize}{0}$\Delta'''$}
\psfrag{gi}{\fontsize{\figurefontsize}{0}$\gamma_i$}
\includegraphics[scale=0.45]{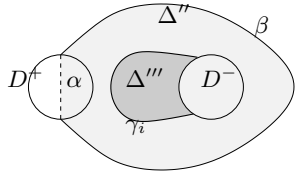} \caption{\small The disc
$\Delta'''$ used to simplify $D\cap D'$.} \label{F:delta3}
\end{figure}

Once again in $Q$, an isotopy may be use to pull $\gamma_i$ and
the whole $\partial D'\cap\Delta'''$ along $\Delta'''$ through
$D$.  This process either
\begin{itemize}

    \item maintains $|\Delta\cap D|$ unchanged removing an arc of
    $\Delta\cap D$ but introducing a closed curve, which we remove
    using the argument previously presented, or

    \item two arcs of $\Delta\cap D$ are joined to one.
\end{itemize}
In both situations we reduce the number of components (all arcs) of $D\cap D'$.

It remains to consider the case when no $\gamma_i$ has both points
in $\partial D^-$. We shall show that this case cannot happen.

Recall that the ends of each $\gamma_i$ have an order induced by the
orientation in $\partial D'$ and note that $\gamma_0=\beta$ has
both ends in $\partial D^+$. It follows that $\gamma_1$ has its
first end in $\partial D^-$. Since no $\gamma_i$ has both ends in
$\partial D^-$ the final end of $\gamma_1$ is in $\partial D^+$.
The same argument shows that if $\gamma_i$ goes from $\partial
D^-$ to $\partial D^+$ then so does $\gamma_{i+1}$, therefore, by
induction, all $\gamma_i$, $i\geq 1$ have this property. But
$\Gamma$ is finite with cyclic order, hence some
$\gamma_n=\gamma_0=\beta$, which does not have this property, a
contradiction.

This completes the analysis of case 2, showing that, also in this
situation, we can choose $D'$ in its isotopy class in such a way
that $D\cap D'=\emptyset$.

The above case analysis shows that we can always
assume that $D'$ is disjoint from $D$. It remains to prove that
$D$ and $D'$ are parallel, hence isotopic. Again we divide the
argument into the same two cases:

1. $F$ is disconnected ($D$ separates $Q$). Let $Q_0$ be the
component of $Q'$ ($Q$ cut open along $D$) containing $D'$. As
before, $Q_0\simeq F_0\times I$ and there is a disc $D_0\subseteq
F_0\times\{1\}$ corresponding to $D$. Since $\partial D'\subseteq
F_0\times\{1\}$, incompressibility of $F_0\times\{1\}$ and
irreducibility of $F_0\times I$ imply that $D'$ is parallel to a
disc $\Delta\subseteq F_0\times\{1\}$. There are only two
possibilities: either $\Delta\cap D_0=\emptyset$ or $\Delta\cap
D_0=D_0$. If $\Delta\cap D_0=\emptyset$ then $D'$ is parallel to
$\partial_e Q$ in $Q$, hence it is not a compressing disc for
$\partial_e Q$, a contradiction. Therefore $\Delta\cap D_0=D_0$,
and $\Delta-D_0$ is an annulus, showing that $D'$ parallel to
$D_0$. Going back to $Q$, $D'$ is parallel to $D$, therefore they
are isotopic.

2. $F$ is connected ($D$ does not separate $Q$). The argument is
analogous to the one in the first case: we consider $Q'\simeq
F\times I$, with two discs $D^+$, $D^-\subseteq Q'$ corresponding
to $D$. Hence $D'\subseteq Q'$ is parallel to $\Delta\subseteq
F\times\{1\}$ and the case $\Delta\cap(D^+\cup D^-)=\emptyset$
cannot happen either. The cases $\Delta\cap(D^+\cup D^-)=D^+$ or
$\Delta\cap(D^+\cup D^-)=D^-$ again give $D'$ parallel to $D$ in
$Q$ (hence isotopic). The only situation not analogous to the
previous case is when $\Delta\cap(D^+\cup D^-)=D^+\cup D^-$. Here,
as can be seen by going back to $Q$, $D'$ separates $Q$, implying
that $F$ is not connected, a contradiction.
\end{proof}

\begin{theorem}\label{T:reduction}
Let $f\colon H\to H$ be a reducible automorphism of a \hndl\ $H$.
If $\partial f=f\vert_{\partial H}$ is pseudo-Anosov, then a closed
reducing surface $F$ satisfies
$$
\chi(\partial H)+4\leq\chi(F)\leq 0.
$$
\end{theorem}

\begin{proof}
Indeed, since $\partial f$ is pseudo-Anosov, $f$ either is \irred\
or, as is the case, admits a closed reducing surface. A closed
reducing surface $F\subseteq\mathring{H}$ is the interior boundary
of a non-trivial \compression\ $(Q,F)$ which is $f$-invariant and
whose exterior boundary is $\partial_eQ=\partial H$. Since
$\partial H$ is connected then so is $\partial_e Q$ and $Q$ also.
By Lemma~\ref{L:euler}
$$
\chi(F)=\chi(\partial H)+2n,
$$
where $n\geq 1$ is the smallest number of $1$-handles of $Q$. Let
$D$ be the disc dual to a $1$-handle. If $n=1$ then
Proposition~\ref{P:1-handle} above implies that $f(D)=D$, hence
$f(\partial D)=\partial D\subseteq
\partial H$, contradicting the hypothesis that $\partial f$ is
pseudo-Anosov. Therefore $n\geq 2$. It is clear that $\chi(F)\leq
0$ for, by definition, spheres are not reducing surfaces.
\end{proof}

\begin{corollary}\label{C:genus2}
Let $f\colon H\to H$ restrict to $\partial H$ as a pseudo-Anosov
automorphism. If the genus of $H$ is two then $f$ is \irred.
\end{corollary}
\begin{proof}
If there were a closed reducing surface $F$ for $f$ then
$2\leq\chi(F)\leq 0$.
\end{proof}

\begin{obs}
This result enables us to reduce the problems of identification or
construction of \irred\ automorphisms of genus two \hndls\ to the
better understood analogues for pseudo-Anosov automorphisms of
surfaces. For instance, Penner's Theorem~\ref{T:penner}
\cite{RP:88} becomes a method for generating \irred\ automorphisms
(e.g. Examples \ref{E:first example} and \ref{E:genus2} below).
Moreover Bestvina and Handel's algorithm to decide whether a given
surface automorphism is pseudo-Anosov or not \cite{BH:Surfaces}
may be used as an algorithm to decide whether an automorphism of a
genus two \hndl\ is \irred\ or not.

On the other hand the same result exposes differences between the
two dimensions (see the remark after Example~\ref{E:genus2}).
\end{obs}

\begin{example}\label{E:genus2}
Consider $H$ a genus 2 \hndl. Figure \ref{F:genus2} shows two
curves $C_0$, $C_1$ in $\partial H$.

\begin{figure}[ht]
\centering \includegraphics[scale=0.25]{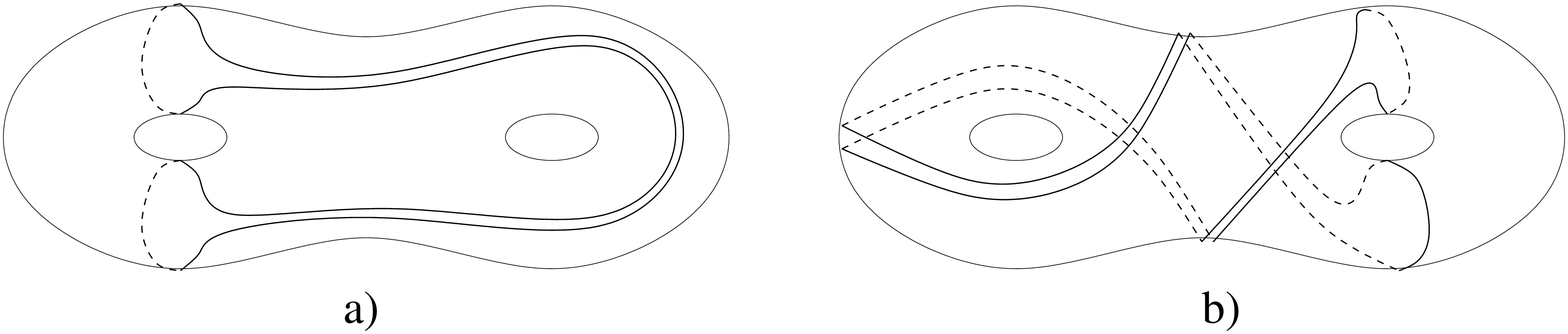}
\caption{\small
        a) the curve $C_0$ bounding the disc $D_0\subseteq H$;
        b) the curve $C_1$ bounding the disc $D_1\subseteq H$.}
\label{F:genus2}
\end{figure}

It is easy to see that these curves bound discs in $H$, $D_0$ and
$D_1$ respectively (one can see them in the picture as the band
sum of discs dual to the handles). We will define $f\colon H\to H$
as the composition of Dehn twists along these discs, to the left
along $D_1$ and to the right along $D_0$:
$$
f=T_{D_0}^+\circ T_{D_1}^-.
$$
It is routine to verify that $(\{C_0\}, \{ C_1\})$ is a Penner
pair in
 $\partial H$, hence
$$
\partial f=T_{C_0}^+\circ T_{C_1}^-
$$
is pseudo-Anosov (Theorem~\ref{T:penner}). By
Corollary~\ref{C:genus2}, $f$ is \irred.

\end{example}

\begin{obs}
The example of \irred\ automorphism $f\colon H\to H$ above is
given as a composition of twists on discs. Since twists on discs
induce the identity on the fundamental group, it is immediate that
$f_*\colon\pi_1(H)\to\pi_1(H)$ is the identity. Therefore, for a
general \irred\ automorphism $f$, $f_*$ may fail to capture its
complexity. We note that that is not the case for pseudo-Anosov
automorphisms of surfaces. This difference should not be regarded
as weakening the analogy between these two classes of
automorphisms, but rather as exposing the richness of the
three-dimensional setting.

\end{obs}

The example below shows that neither inequalities of Theorem
\ref{T:reduction} can be improved.

\begin{example}\label{E:genus3}
As in Example~\ref{E:genus2}, Figure~\ref{F:genus3} a), b)
represent the boundaries of two discs in a handlebody $H$, here
with genus 3. These boundaries yield a Penner pair in $\partial H$
as well. Hence a composition $f$ of twists to opposite directions
along these discs yields $\partial f$ as a pseudo-Anosov
automorphism. Theorem~\ref{T:reduction} says that the Euler
characteristic of a closed reducing surface is zero. Indeed, one
can see that there exists a torus that does not intersect the
discs -- therefore being invariant under $f$
(Figure~\ref{F:genus3} c)).

\begin{figure}[h]
\centering
\includegraphics[scale=0.18]{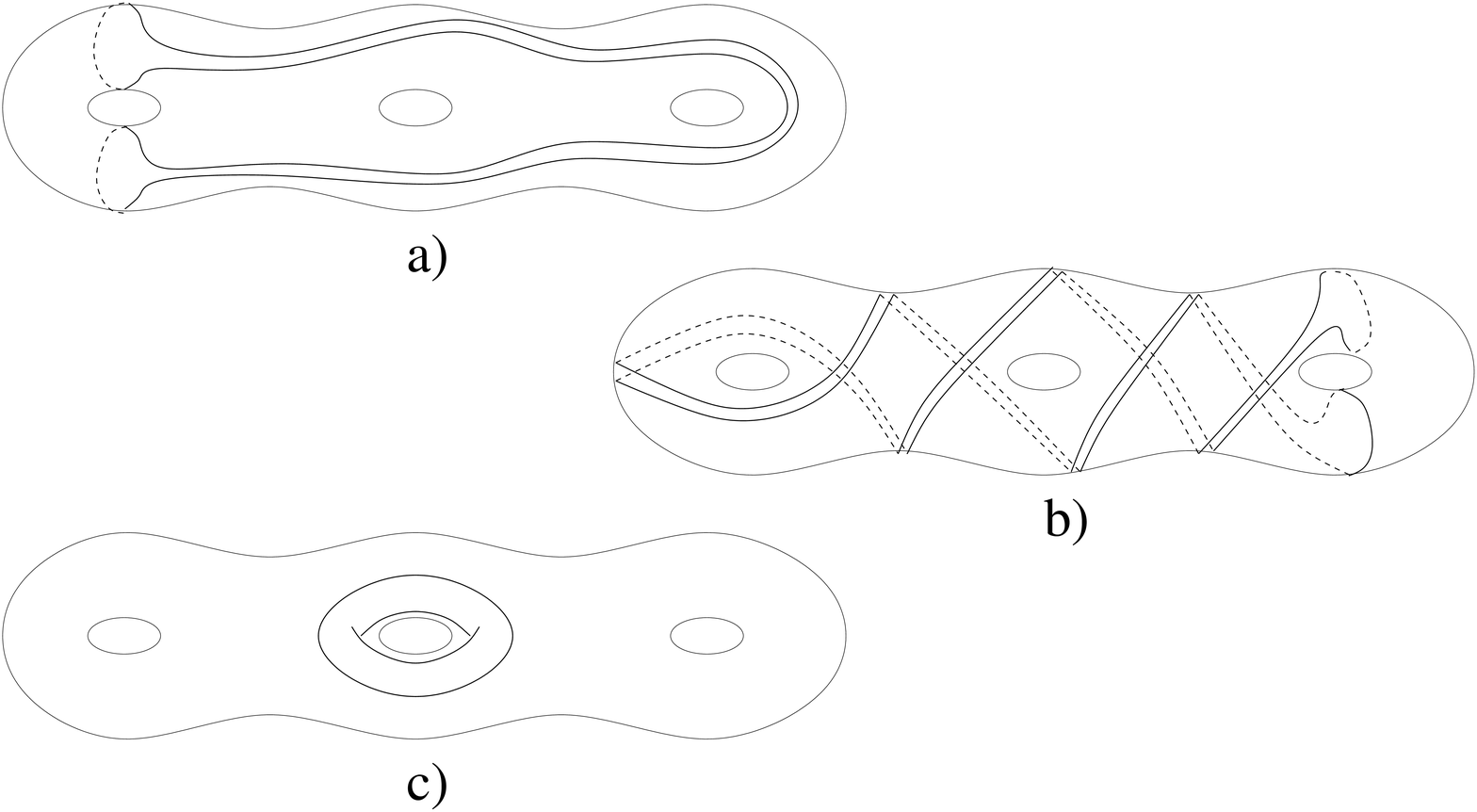}
\caption{\small
        a) the curve $C_0$ bounding a disc $D_0\subseteq H$;
        b) the curve $C_1$ bounding a disc $D_1\subseteq H$;
        c) the invariant torus.}
\label{F:genus3}
\end{figure}

\end{example}

\bibliographystyle{halpha}
\bibliography{reference}

\end{document}